\documentclass[12pt]{article}  
\def\sq{\hbox {\rlap{$\sqcap$}$\sqcup$}}
\def\sq{\hbox {\rlap{$\sqcap$}$\sqcup$}}
\def\R{ {\rm R \kern -.31cm I \kern .15cm}}
\def\C{ {\rm C \kern -.15cm \vrule width.5pt \kern .12cm}}
\def\Z{ {\rm Z \kern -.27cm \angle \kern .02cm}}
\def\N{ {\rm N \kern -.26cm \vrule width.4pt \kern .10cm}}
\def\1{{\rm 1\mskip-4.5mu l} }
\def\lsim{\raise0.3ex\hbox{$<$\kern-0.75em\raise-1.1ex\hbox{$\sim$}}}
\def\gsim{\raise0.3ex\hbox{$>$\kern-0.75em\raise-1.1ex\hbox{$\sim$}}}
\def\noi{\noindent}

\def\beq{\begin{equation}}   \def\eeq{\end{equation}}
\def\bea{\begin{eqnarray}}  \def\eea{\end{eqnarray}}
\def\nn{\nonumber}
\def\noi{\noindent}
\def\beeq{\begin{eqnarray}} \def\eeeq{\end{eqnarray}}
\newcommand\mysection{\setcounter{equation}{0}\section}

\newcounter{hran}

\begin{document} 
\centerline{\large\bf Modified wave operators without loss of regularity} 
 \vskip 3 truemm 
 \centerline{\large\bf  for some long range Hartree equations. II} 
  \vskip 0.8 truecm

\centerline{\bf J. Ginibre}
\centerline{Laboratoire de Physique Th\'eorique\footnote{Unit\'e Mixte de
Recherche (CNRS) UMR 8627}}  \centerline{Universit\'e de Paris XI, B\^atiment
210, F-91405 Orsay Cedex, France}
\vskip 3 truemm

\centerline{\bf G. Velo}
\centerline{Dipartimento di Fisica e Astronomia, Universit\`a di Bologna}  \centerline{and INFN, Sezione di
Bologna, Italy}

\vskip 1 truecm

\begin{abstract}
We continue the study of the theory of scattering for some long range Hartree equations with potential $|x|^{-\gamma}$, performed in a previous paper, denoted as I, in the range $1/2 < \gamma < 1$. Here we extend the results to the range $1/3 < \gamma < 1/2$. More precisely we study the local Cauchy problem with infinite initial time, which is the main step in the construction of the modified wave operators. We solve that problem without loss of regularity between the asymptotic state and the solution, as in I, but in contrast to I, we are no longer able to cover the entire subcriticality range of regularity of the solutions. The method is an extension of that of I, using a better approximate asymptotic form of the solutions obtained as the next step of a natural procedure of successive approximations.
\end{abstract}

\vskip 1 truecm
\noi 2000 MSC :  Primary 35P25. Secondary 35B40, 35Q40, 81U99.\par \vskip 2 truemm

\noi Key words : Long range scattering, wave operators, Hartree equation. \par 
\vskip 1 truecm

\noindent LPT Orsay 13-16\par
\noindent February 2013\par \vskip 3 truemm

\newpage
\pagestyle{plain}
\baselineskip 18pt
\mysection{Introduction}
\hspace*{\parindent} 
This paper is the continuation of a previous paper \cite{5r}, hereafter referred to as I, where we studied the theory of scattering and more precisely the proof of existence of modified wave operators for the long range Hartree type equation 
\beq
i \partial_t u = - (1/2) \Delta u + g(u) u \ . 
\label{1.1e}
\eeq

\noi Here $u$ is a complex valued function defined in space time ${I\hskip-1truemm R}^{n+1}$ with $n \geq 2$, 
$\Delta$ is the Laplace operator in ${I\hskip-1truemm R}^n$ and 
\beq
\label{1.2e}
g(u) = \kappa \  |x|^{-\gamma} \star |u|^2
\eeq

\noi where $\kappa \in {I\hskip-1truemm R}$, $0 < \gamma \leq 1$ and $\star$ denotes the convolution in ${I\hskip-1truemm R}^n$. \par

The main step of that existence proof consists in solving the local Cauchy problem with infinite initial time for (\ref{1.1e}), namely to construct solutions thereof with prescribed asymptotic behaviour as $t \to \pm \infty$. We refer to the introduction of I and to \cite{2r} \cite{3r} \cite{6r} \cite{7r} for general background. In the long range situation $\gamma \leq 1$ that we consider, the asymptotic behaviour of $u$ includes a phase which diverges at infinity in time, and is parametrized by an asymptotic state which plays the role of initial data at infinity. In \cite{2r}, we solved the previous local Cauchy problem at infinity in the range $1/2 < \gamma < 1$ (the easier borderline case $\gamma = 1$ can be treated by the same method), but the treatment in \cite{2r} involved a loss of regularity between the asymptotic state and the solution and failed to cover the entire natural subcritical range of regularity for the asymptotic state. These two defects were remedied in \cite{6r} and \cite{7r} in the cases $\gamma = 1$ and $1/2 < \gamma < 1$ respectively. The main results of \cite{6r} \cite{7r} were then recovered in $I$ by what we regard as a simpler method combining ingredients of \cite{2r} and \cite{6r}. On the other hand, the problem was solved in \cite{3r} by an extension of the method of \cite{2r} for $\gamma \leq 1/2$, actually for any $\gamma \leq 1$, again with a loss of regularity between the asymptotic state and the solution. That loss increases as $\gamma$ decreases through inverse integer values. Now it turns out that the simple method of I can be extended below $\gamma = 1/2$ to solve the problem without loss of regularity. For $\gamma < 1/2$ however, it no longer allows to cover the whole subcritical range, and stronger regularity of the asymptotic state is needed. Furthermore the treatment, although still elementary, becomes increasingly cumbersome as $\gamma$ decreases through inverse integer values. The present paper is devoted to the application of that method to the next accessible range, namely for $1/3 < \gamma < 1/2$, as an illustration of that possibility. The easier borderline case $\gamma = 1/2$ can be treated by the same method, but requires slightly different formulas.\par

We now introduce the relevant parametrization of $u$ needed to study the Cauchy problem at infinite time, restricting our attention to positive time. The unitary group 
\beq
\label{1.3e}
U(t) = \exp (i(t/2)\Delta ) \eeq

\noi which solves the free Schr\"odinger equation can be written as 
\beq
U(t) = M(t) \ D(t) \ F\ M(t) 
\label{1.4e}
\eeq

\noi where $M(t)$ is the operator of multiplication by the function
\beq
M(t) = \exp (ix^2/2t) \ , 
\label{1.5e}
\eeq

\noi $F$ is the Fourier transform and $D(t)$ is the dilation operator
\beq
D(t) = (it)^{-n/2} \ D_0(t) 
\label{1.6e}
\eeq

\noi where
\beq
\left ( D_0 (t) f \right ) (x) = f(x/t)\ . 
\label{1.7e}
\eeq

\noi For any function $w$ of space time, we define
\beq
\label{1.8e}
\widetilde{w}(t) = U(-t) \ w(t) 
\eeq

\noi and we define the pseudoconformal inverse $w_c$ of $w$ by
\beq
w(t) = M(t) \ D(t) \ \overline{w_c} (1/t) 
\label{1.9e}
\eeq

\noi or equivalently
\beq
\widetilde{w}(t) = \overline{F \widetilde{w}_c (1/t)}\ 
\label{1.10e}
 \eeq

\noi which shows that the pseudoconformal inversion is involutive. \par

The Cauchy problem at infinite initial time for $u$ is then equivalent to the Cauchy problem at initial time zero for its pseudoconformal inverse $u_c$. The equation (\ref{1.1e}) is replaced by 
\beq
i \partial_t  u_c = - (1/2) \Delta u_c + t^{\gamma -2} \ g(u_c) u_c \ .
\label{1.11e}
 \eeq

\noi We now parametrize $u_c$ in terms of an amplitude $v$ and a phase $\varphi$ according to 
\beq
\label{1.12e}
u_c(t) = \exp (- i \varphi (t)) v(t) 
\eeq

\noi so that
\begin{eqnarray*}
u(t) &=& M(t)\ D(t) \exp \left ( i \varphi (1/t)\right ) \overline{v}(1/t) \\
&=& D(t) \exp \left ( i \varphi (1/t)\right ) D^{-1}(t)\ M(t)\ D(t) \ \overline{v}(1/t)
\end{eqnarray*}

\noi or equivalently
\beq
\label{1.13e}
u(t) = \exp \left ( i \left ( D_0(t) \varphi (1/t)\right ) \right ) v_c(t)\ .
\eeq

\noi The original equation then becomes the following equation for $v$
\beq
\label{1.14e}
i \partial_t v = - (1/2) \Delta_s v + \left ( t^{\gamma -2} g(v) - \partial_t \varphi \right ) v 
\eeq

\noi where $s = \nabla \varphi$ and
\beq
\Delta_s = (\nabla - i\ s)^2 = \Delta - 2 i\ s \cdot \nabla - i (\nabla \cdot s) - |s|^2 \ .
\label{1.15e}
\eeq

We want to choose $\varphi$ so as to cancel the divergence at $t = 0$ of the last term in (\ref{1.14e}), but that cancellation is needed only at large distances, namely for low momentum. We therefore introduce a momentum cut-off as follows. Let $\chi \in {\cal C}^\infty ({I\hskip-1 truemm R}^+,{I\hskip-1 truemm R}^+)$, $0 \leq \chi \leq 1$, $\chi (\ell ) = 1$ for $\ell \leq 1$, $\chi (\ell ) = 0$ for $\ell \geq 2$. We define 
\beq
\label{1.16e}
\chi_L = \chi (\omega t^{1/2}) \quad , \quad \chi_S = 1 - \chi_L
\eeq

\noi with $\omega = (-\Delta )^{1/2}$, and correspondingly 
\beq
\label{1.17e}
g_L(v) = \chi_L \ g(v) \quad , \quad g_S(v) = \chi_S \ g(v) \ .
\eeq

\noi We want to solve (\ref{1.14e}) with $v$ continuous at $t=0$ with $v(0) = v_0$ for a given $v_0$. In I we chose $\varphi$ such that 
\beq
\label{1.18e}
\partial_t \varphi -  t^{\gamma - 2} \ g_L (v_0) \approx 0\ ,
\eeq

\noi a choice which was sufficient for $\gamma > 1/2$. That choice however is not sufficient for $\gamma \leq 1/2$ since then the terms coming from $|s|^2$ and from $g_L(v) - g_L (v_0)$ in (\ref{1.14e}) both fail to be integrable at $t=0$. We must therefore choose a better asymptotic form $v_a$ for $v$, still with $v_a (0) = v_0$. We rewrite (\ref{1.14e}) as 
\beq
\label{1.19e}
i \partial_t v = L(v)v
\eeq

\noi with
\bea
\label{1.20e}
L(v) &=& - (1/2) \Delta + is\cdot \nabla + (i/2)(\nabla \cdot s) + t^{\gamma -2} g_S(v) \nn \\
&+& t^{\gamma -2} \left ( g_L (v) - g_L (v_a)\right ) + (1/2) |s|^2 + t^{\gamma -2} g_L(v_a) - \partial_t \varphi \ .
\eea

\noi If the asymptotic $v_a$ is sufficiently accurate, we may expect that the term with $g_L(v) - g_L(v_a)$ will be integrable at $t=0$ and we may try to cancel the remaining divergences by choosing $\varphi$ according to 
\beq
\label{1.21e}
\partial_t \varphi = t^{\gamma - 2} g_L (v_a)  + (1/2) |s|^2\ ,  
\eeq

\noi with initial condition $\varphi (1) = 0$, since the RHS of (\ref{1.21e}) fails to be integrable at $t = 0$. In order to control the term with $g_L (v) - g_L(v_a)$, as in I and following \cite{6r}, we use the facts that $g(v)$ depends only on $|v|^2$ and that, if $v$ satisfies a linear Schr\"odinger equation 
\beq
\label{1.22e}
i \partial_t v = - (1/2) \Delta_s v + V\ v   
\eeq

\noi for some real potential $V$, then $v$ satisfies the local conservation law
\beq
\label{1.23e}
\partial_t |v|^2 = - \ {\rm Im} \ \overline{v} \Delta v + \nabla \cdot s |v|^2 \ .  
\eeq

\noi If we impose that $v_a$ satisfies the transport equation
\beq
\label{1.24e}
\partial_t v_a = s\cdot \nabla v_a + (1/2) (\nabla\cdot s) v_a \ ,  
\eeq

\noi then we obtain
\beq
\label{1.25e}
\partial_t \left ( |v|^2  - |v_a|^2\right ) = -  {\rm Im} \ \overline{v} \Delta v + \nabla \cdot s \left ( |v|^2  - |v_a|^2\right ) 
\eeq

\noi which provides a good starting point to estimate $g_L(v) - g_L (v_a)$. (One could also impose the Schr\"odinger equation 
$$i \partial_t v_a = - (1/2) \Delta_s v_a$$

\noi but that would introduce unnecessary complications without improving the crucial estimates).\par

We are therefore led to choose $(\varphi , v_a)$ by solving the system (\ref{1.21e}) (\ref{1.24e}) with initial conditions $\varphi (1) = 0$, $v_a(0) = v_0$, but this is a nonlinear system which is hardly simpler than the original equation, and we seem to have gained nothing so far. However $(\varphi , v_a)$ are only asymptotic quantities, and it suffices to solve that system approximately by iteration. We therefore define successive approximate solutions $(\varphi_m , v_m)$ of (\ref{1.21e}) (\ref{1.24e}) by
\beq
\label{1.26e}
 \left\{ \begin{array}{l}
\partial_t \varphi_m = t^{\gamma - 2} g_L (v_m) + (1/2) |s_{m-1}|^2 \\
\\
\partial_t v_m = s_{m-1} \cdot \nabla v_m + (1/2) \left ( \nabla \cdot s_{m-1}\right ) v_m 
\end{array} \right .
\eeq

\noi with $\varphi_m(1) = 0$, $v_m(0) = v_0$. The system (\ref{1.26e}) determines $v_m$ by a linear transport equation with a smooth vector field $s_{m-1}$  and then $\varphi_m$ by integration over time. The choice made in I was essentially the case $m = 0$ (with $\varphi_{-1} \equiv 0$), namely
\beq
\label{1.27e}
 \left\{ \begin{array}{l}
\partial_t \varphi_0 = t^{\gamma - 2} g_L (v_0) \\
\\
\partial_t v_0 = 0 
\end{array} \right .
\eeq

\noi where by a slight abuse of notation we denote by $v_0$ both the initial value $v_m (0)$ and the constant function of time equal to $v_0$. That choice was adequate for $\gamma > 1/2$. In the present paper, we use the next approximation $m = 1$, namely we take $(\varphi , v_a) = (\varphi_1, v_1)$ so that 
$$
 \hskip 4 truecm \left\{ 
\begin{array}{ll}
\partial_t \varphi = t^{\gamma - 2} g_L (v_a) + (1/2) |s_0|^2 \hskip 3.9 truecm &{\rm (1.28)} \\
\\
\partial_t v_a = s_0 \cdot \nabla v_a + (1/2) (\nabla \cdot s_0) v_a &{\rm (1.29)}
\end{array} 
\right .
$$

\noi with $s_0 = \nabla \varphi_0$ defined by (\ref{1.27e}), with $\varphi_0 (1) = \varphi (1) = 0$ and $v_a(0) = v_0$. That choice turns out to be sufficient to cover the range $1/3 < \gamma \leq 1/2$. With that choice, the basic equation to be solved is (\ref{1.19e}), where now
$$L(v) = - (1/2) \Delta + i s \cdot \nabla + (i/2) (\nabla \cdot s ) + t^{\gamma - 2} g_S (v)$$
$$+t^{\gamma - 2} \left ( g_L (v) - g_L (v_a)\right ) + (1/2) \left ( |s|^2 - |s_0|^2 \right )\eqno(1.30)$$

\noi and 
$$s = s_0 + s_b + s_c \ , \eqno(1.31)$$
$$s_0 = - \nabla \int_t^1 dt' \ t{'}^{\gamma - 2} g_L (v_0) \ , \eqno(1.32)$$
$$s_b = - (1/2) \nabla  \int_t^1 dt' |s_0 (t')|^2\ , \eqno(1.33)$$
$$s_c = - \nabla  \int_t^1 dt' \ t{'}^{\gamma - 2} \left ( g_L (v_a) - g_L (v_0)\right ) (t') \ . \eqno(1.34)$$

\noi More generally, one expects the approximation $(\varphi_m, v_m)$ to be sufficient to cover the range $1/(m+2) < \gamma \leq 1/(m+1)$.\par

In the present paper, we treat the problem in the range $1/3 < \gamma < 1/2$ with the previous choice of $(\varphi , v_a)$. (The simpler case $\gamma = 1/2$ can be treated with the same choice, but requires slightly different formulas). As mentioned above, we solve the local Cauchy problem at infinity in time for $u$ (at time zero for $v$) without any regularity loss between the asymptotic state $v_0$ and the solution $v$, but in contrast to I, we are unable to cover the entire subcritical range for $v$, and stronger than subcritical regularity is required for $v_0$ as soon as $\gamma < 1/2$. \par

In addition to (\ref{1.19e}), we shall also need the partly linearized equation for $v'$
$$i \partial_t v' = L(v)v' \eqno(1.35)$$

\noi with $L(v)$ again defined by (1.30).\par

The method consists in first solving the Cauchy problem with initial time zero for the linearized equation (1.35). One then shows that the map $v \to v'$ thereby defined is a contraction in a suitable space in a sufficiently small time interval. This solves the Cauchy problem with initial time zero for the nonlinear equation (1.19). One then translates the results through the change of variables (1.12) to solve the Cauchy problem with initial time zero for the equation (1.11) or equivalently with infinite initial time for the equation (1.1). The final result can be stated as the following proposition, which is adapted to the equation (\ref{1.1e}) in a neighborhood of infinity in time. We need the notation
$$FH^\rho = \{ u \in {\cal S} ' : F^{-1} u \in H^\rho \} \ .$$

\noi {\bf Proposition 1.1.} {\it Let $1/3 < \gamma < 1/2$ and $2-5\gamma /2 < \rho < n/2$.  Let $u_0 \in F H^\rho$. Let $\varphi$ be defined by (1.28) with $\varphi (1) = 0$, with $s_0$ defined by (1.32) and $v_a$ defined by (1.29) with $v_a (0) = v_0$ (see Lemma 2.5, part (2)). Then there exists $T_\infty > 0$ and there exists a unique solution $u$ of the equation (\ref{1.1e}) such that $v_c$ defined by (\ref{1.9e}) (\ref{1.12e}) or equivalently by (\ref{1.13e}) satisfies $\widetilde{v}_c \in {\cal C} ([T_\infty , \infty ), FH^\rho )$ and such that
$$
\widetilde{v}_c (t) \to u_0 \hbox{ in } FH^\rho \hbox{ when } t \to \infty \ . \eqno(1.36)
$$

\noi Furthermore $\widetilde{u} \in {\cal C} ([T_\infty , \infty ), FH^\rho )$ and $\widetilde{u}$ satisfies the estimate 
$$\parallel \widetilde{u} (t) ; FH^{\rho} \parallel \ \leq C\ a_0 \left ( 1 + a_0^2 (1 + a_0^2)  t^{1 - \gamma} \right )^{1 + [\rho ]} 
\eqno(1.37)$$

\noi for all $t \geq T_\infty$, where $[\rho ]$ is the integral part of $\rho$ and} 
$$a_0 = \ \parallel u_0 ; FH^\rho \parallel \ .$$
\vskip 3 truemm

Proposition 1.1 follows from Propositions 4.1 and  4.2 through the change of variables (\ref{1.9e}) or (\ref{1.10e}), which implies in particular that
$$\parallel \widetilde{w} (t) ; FH^\rho \parallel \ = \ \parallel \widetilde{w}_c (1/t) ; H^\rho \parallel \ = \ \parallel w_c (1/t) ; H^\rho \parallel \ .
\eqno(1.38)$$

\noi As previously mentioned, the condition $\rho > 2 - 5\gamma /2 = 1 - \gamma /2 + 1 - 2 \gamma$ is stronger thant the subcriticality condition $\rho > 1 - \gamma /2$ for $\gamma < 1/2$. The technical origin of that condition is explained in Remark 3.3 below. \par

This paper is organized as follows. In Section~2, we introduce some notation and we collect a number of estimates which are used throughout this paper. In Section~3, we study the Cauchy problem for the linearized equation (1.35) with initial time $t_0 \geq 0$. In Section~4, we solve the Cauchy problem with initial time zero for the nonlinear equation (\ref{1.19e}) and we translate the result into the corresponding one for the equation (1.11).\par

This paper follows I closely and uses the same methods. In order to make it reasonably self contained while avoiding excessive repetition, we have given full statements of the intermediate results, but we have shortened or even omitted some of the proofs when they are identical with those of I.

\mysection{Notation and preliminary estimates} 
\hspace*{\parindent}
In this section we introduce some notation and we collect a number of estimates which will be used throughout this paper. We denote by
$\parallel \cdot \parallel_r$ the norm in $L^r \equiv L^r({I\hskip-1truemm R}^{n})$.  For any interval $I$ and any Banach space $X$ we denote by ${\cal C}(I,X)$ (resp. ${\cal C}_w(I,X))$ the space of strongly (resp. weakly) continuous functions from $I$ to $X$ and by $L^{\infty} (I, X)$ the space of measurable essentially bounded functions from $I$ to $X$. For real numbers $a$ and $b$ we use the notation $a \vee b = {\rm Max}(a,b)$ and $a\wedge b = {\rm Min} (a,b)$. We define $(a)_+ = a \vee 0$ and 
$$\begin{array}{lll}
[a]_+ &= (a)_+ &\qquad \hbox{for $a\not= 0$} \\
&& \\
&= \varepsilon \hbox{ for some $\varepsilon > 0$} &\qquad \hbox{for $a= 0$} \ .
\end{array}$$
\par

We shall use the Sobolev spaces $\dot{H}_r^\sigma$ and $H_r^\sigma$ defined for $- \infty < \sigma < + \infty$, $1 \leq r < \infty$ by
$$\dot{H}_r^\sigma = \left \{ u:\parallel u;\dot{H}_r^\sigma\parallel \ \equiv \ \parallel \omega^\sigma u\parallel_r \ <
\infty \right \}$$

\noi and

$$H_r^\sigma = \left \{ u:\parallel u;H_r^\sigma\parallel \ \equiv \ \parallel <\omega>^\sigma u\parallel_r \ <
\infty \right \}$$

\noi where $\omega = (- \Delta)^{1/2}$ and $< \cdot > = (1 + |\cdot |^2)^{1/2}$. The subscript $r$ will
be omitted both in $H^\sigma$ and in the $L^r$ norm if $r = 2$ and we shall use
the notation 
$$\parallel \omega^{\sigma \pm 0} u \parallel \ = \left ( \parallel \omega^{\sigma + \varepsilon} u \parallel\ \parallel \omega^{\sigma - \varepsilon} u \parallel\right )^{1/2} \quad \hbox{for some $\varepsilon > 0$\ .}$$  

Note also that for $0 < \gamma < n$ \cite{8r} 
$$g(u) = \kappa \ |x|^{-\gamma} \star |u|^2 = \kappa \ C_{\gamma , n} \ \omega^{\gamma - n} \ |u|^2 \ .$$ 

We shall use extensively the following Sobolev inequalities. \\

\noi {\bf Lemma 2.1.} {\it Let $1 < q, r < \infty$, $1 < p \leq \infty$ and $0 \leq \sigma < \rho$. If $p = \infty$, assume that $\rho - \sigma > n/r$. Let $\theta$ satisfy $\sigma /\rho \leq \theta \leq 1$ and
$$n/p - \sigma = (1 - \theta )n/q + \theta (n/r - \rho ) \ .$$

\noi Then the following inequality holds}
\beq  
\label{2.1e}
\parallel \omega^\sigma u\parallel_p \ \leq C  \parallel u \parallel_q^{1 - \theta} 
 \ \parallel \omega^\rho u \parallel_r^\theta \ .
 \eeq

We shall also use extensively the following Leibnitz estimates. \\

\noi {\bf Lemma 2.2.} {\it Let $1 < r,r_1,r_3 < \infty$ and
$$1/r = 1/r_1 + 1/r_2 = 1/r_3 + 1/r_4 \ .$$

\noi Then the following estimates hold for $\sigma \geq 0$~:}
\beq  
\label{2.2e}
\parallel \omega^\sigma (uv) \parallel_r \ \leq C \left ( \parallel  \omega^\sigma u \parallel_{r_1} 
 \ \parallel v \parallel_{r_2} + \parallel  \omega^\sigma v \parallel_{r_3} \  
 \parallel u \parallel_{r_4} \right )
\ .  \eeq 

An easy consequence of Lemmas 2.1 and 2.2 is the inequality 
\bea  
\label{2.3e}
\parallel \omega^\sigma f\ u \parallel \ &\leq& C  \left ( \parallel f \parallel_\infty \ + \ \parallel \omega^{n/2} f \parallel \right ) 
 \ \parallel  \omega^{\sigma} u \parallel  \  .\nn \\
&\leq& C \parallel \omega^{n/2\pm 0} f \parallel  \ \parallel  \omega^{\sigma} u \parallel 
 \eea

\noi which holds for $|\sigma | < n/2$. \par

Another consequence is the following lemma.\\

\noi {\bf Lemma 2.3.} {\it Let $0 < \sigma = \sigma_1  + \sigma_2$ and $\sigma_1 \vee \sigma_2 < n/2$. Then} 
\beq  
\label{2.4e}
\parallel \omega^{\sigma - n/2} (uv) \parallel \ \leq C  \parallel  \omega^{\sigma_1} u \parallel \ \parallel  \omega^{\sigma_2} v \parallel  \  .
 \eeq

We shall also need some commutator estimates, which are most conveniently stated in terms of homogeneous Besov spaces $\dot{B}_{r, q}^\sigma$ \cite{1r}. In the applications, we shall use only the fact that $\dot{B}_{2,2}^\sigma$ = $\dot{H}^\sigma$. \\

\noi {\bf Lemma 2.4.} {\it Let $P_i$, $i = 1,2$ be homogeneous derivative polynomials of degree $\alpha_i$ or $\omega^{\alpha_i}$ for $\alpha_i \geq 0$. Let $\lambda > 0$. Then for any (sufficiently regular) functions $m$, $u$ and $v$ the following estimates hold.
\bea  
\label{2.5e}
&&|< P_1 u, [\omega^\lambda , m] P_2 v>| \leq  C \parallel  m ; \dot{B}_{r_0,2}^{\sigma_0} \cap \nabla^{-1} \omega^{1-\nu} L^{q_0} \parallel \ \parallel u; \dot{B}_{r_1,2}^{\sigma_1} \cap L^{q_1} \parallel \nn \\
&& \parallel v;  \dot{B}_{r_2,2}^{\sigma_2} \cap L^{q_2}\parallel 
 \eea
\noi with $0 \leq \nu \leq 1$, $1 \leq r_i, q_i \leq \infty$, $0 \leq i \leq 2$,
\beq  
\label{2.6e}
\delta (q_0) = \sigma_0 + \delta (r_0) - \nu \quad , \quad \delta (q_i) = \sigma_i + \delta (r_i) \ , \ i = 1,2 \ .
 \eeq
\beq  
\label{2.7e}
\sum_{0 \leq i \leq 2} \sigma_i + \delta (r_i) = \lambda + \alpha_1 + \alpha_2 + n/2
 \eeq
\beq  
\label{2.8e}
\left\{ \begin{array}{l} \sigma_0  + \left ( \sigma_1  \wedge \sigma_2 \right ) \geq \lambda + \alpha_1 + \alpha_2\\ \\ \sigma_1  +  \sigma_2  \geq \lambda + \alpha_1 + \alpha_2 - \nu \end{array} \right .
 \eeq
\noi where $\delta (r) \equiv n/2 - n/r$ and $\nabla^{-1} \omega^{1-\nu} L^q$ is the space of tempered distributions $m$ such that $\omega^{\nu - 1} \nabla m \in L^q$.}\\

\noi The proof is given in Appendix A1 of I. \par

We shall need some estimates of $s_0$, $v_a$ and $s$ defined by (1.29) and (1.31)-(1.34). For $0 < \gamma < 1$, $0 < \rho < n/2$ and $\alpha \in {I \hskip - 1 truemm R}$, we define 
\beq  
\label{2.9e}
\lambda_\alpha = \gamma - (1/2) [\alpha + 1 + \gamma - 2 \rho ]_+  
 \eeq

\noi so that $\lambda_\alpha$ is a decreasing function of $\alpha$ and $\lambda_\alpha \leq \gamma$ in all cases. The subcriticality condition $\rho > 1 - \gamma /2$ is equivalent to $\lambda_1 > 0$. Furthermore, under that condition, $\lambda_0 = \gamma$ for $\gamma < 1/2$ since then $(1 + \gamma )/2 < 1 - \gamma /2$. For clarity, we shall nevertheless keep $\lambda_0$ explicitly in that case in some of the estimates.\\

\noi {\bf Lemma 2.5.} {\it Let $0 < \gamma < 1$ and $0 < \rho < n/2$. Let $v_0 \in H^\rho$. \par

(1) Let $s_0$ be defined by (1.32). Then
\beq
\label{2.10e}
\parallel \omega^{\alpha + n/2\pm 0}  s_0 \parallel \ \leq \ C  \ a_0^2 \ t^{\lambda_\alpha - 1} 
\eeq

\noi for $\alpha + 1 + \gamma > 0$, where $a_0 = \parallel v_0 ; H^\rho \parallel$. \par

Let in addition $\rho > 1 - \gamma /2$. Let $0 < T \leq 1$, $I = (0, T]$ and $\overline{I} = [0, T]$. Then \par

(2) The equation (1.29) for $v_a$ with initial condition $v_a(0) = v_0$ has a unique solution $v_a \in {\cal C}(\overline{I}, H^\rho )$ and that solution satisfies the estimate
\beq
\label{2.11e}
\parallel v_a ; L^\infty (I, H^\rho )\parallel  \ \equiv \ a_a \leq a_0 \exp \left ( C\ a_0^2 \ T^{\lambda_1}   \right ) \ .
\eeq

\noi Let in addition $\gamma < 1/2$. Then \par

(3) The following estimates hold for $s_b$ and $s_c$ defined by (1.33) (1.34)
\bea
\label{2.12e}
&&\parallel \omega^{\alpha +n/2 \pm 0}  \ s_b \parallel\ \leq \ C \ a_0^4 (1 - 2 \gamma )^{-1} \ t^{\lambda_0 + \lambda_{\alpha + 1} - 1}\\
&&\parallel \omega^{\alpha +n/2 \pm 0}  \ s_c \parallel\ \leq \ C \ a_0^2 \ a_a^2 \gamma^{-1}  (1 - 2 \gamma )^{-1} \ t^{\lambda_0 + \lambda_{\alpha + 1} - 1}
\label{2.13e}
\eea

\noi for $\alpha + 2 + \gamma > 0$ and for all $t \in I$.} \\

\noi {\bf Sketch of proof.}

\noi \underline{Part (1)} follows from the fact that  
\bea
\label{2.14e}
\parallel \omega^{\alpha +n/2 \pm 0}  \ s_0 \parallel &\leq&  C \int_t^1 dt'\ t{'}^{\gamma -2} \parallel \omega^{\alpha + 1 + \gamma - n/2 \pm 0} \ \chi_L(t') |v_0|^2\parallel \nn \\
&\leq&  C \ a_0^2 \int_t^1 dt' \  t{'}^{\lambda_\alpha - 2} 
\eea

\noi by (\ref{2.4e}). \\

\noi \underline{Part (2)}. The existence of a unique solution $v_a$ of (1.29) as stated can be proved easily, for instance by first solving the Cauchy problem with initial condition $v_a (t_0) = v_0$ for some $t_0 > 0$ by a parabolic regularisation, a fixed point argument and a limiting procedure, and then taking the limit $t_0 \to 0$ of the solution thereby obtained. The key technical fact consists of preliminary versions of the a priori estimate (\ref{2.11e}), which we now derive. From (1.29), we obtain 
\bea
\label{2.15e}
&&\partial_t \parallel \omega^{\sigma}  v_a \parallel^2 \ =  \ 2 {\rm Re} \ \langle \omega^\sigma  v_a , \omega^\sigma \left ( s_0 \cdot \nabla v_a + (1/2) (\nabla \cdot s_0) v_a \right ) \rangle \nn \\ 
&&\leq \ C \parallel \omega^{1 + n/2 \pm 0} s_0(t) \parallel  \ \parallel  \omega^{\sigma} v_a \parallel^2 \nn \\
&&\leq \ C \ a_0^2 \ t^{\lambda_1 - 1} \parallel \omega^{\sigma}  v_a \parallel^2
\eea

\noi for $0 < \sigma \leq \rho$, by Lemma 2.4 and (\ref{2.10e}), from which (\ref{2.11e}) follows by integration over time.\\

\noi \underline{Part (3)}. We first estimate $s_b$. From (1.33) we obtain
\bea
\label{2.16e}
&&\parallel \omega^{\alpha +n/2 \pm 0}  \ s_b \parallel \ \leq\   \int_t^1 dt' \parallel \omega^{\alpha + 1 + n/2 \pm 0} \ |s_0 (t')|^2\parallel \nn \\
&&\leq \ \int_t^1 dt'   \parallel \omega^{\alpha + 1 + n/2 \pm 0} \ s_0 \parallel \ \parallel \omega^{n/2 \pm 0}  \ s_0 \parallel (t')
\eea

\noi by (\ref{2.3e}), for $\alpha + 1 + n > 0$,
$$\leq \ C\ a_0^4 \int^1_t dt' \ t{'}^{\lambda_0 + \lambda_{\alpha + 1} - 2}$$

\noi by (\ref{2.10e}), for $\alpha + 2 + \gamma > 0$,
\beq
\label{2.17e}
\leq \ C \ a_0^4 \left ( 1 - \lambda_0 -  \lambda_{\alpha + 1} \right )^{-1}   \ t^{\lambda_0 + \lambda_{\alpha + 1} - 1} 
\eeq

\noi from which (\ref{2.12e}) follows since $\lambda_0$, $\lambda_{\alpha + 1} \leq \gamma$. The condition $\alpha + 2 + \gamma > 0$ implies $\alpha + 1 + n > 0$ since $\gamma \leq 1$.\par

We next estimate $s_c$. From (1.34) and from the conservation law
\beq
\label{2.18e}
\partial_t |v_a|^2 = \nabla \cdot s_0 |v_a|^2
\eeq

\noi with $v_a (0) = v_0$, we obtain
\beq
\label{2.19e}
s_c = - \int_t^1 dt' \ t{'}^{\gamma -2} \int_0 ^{t'} dt'' \ \nabla \omega^{\gamma - n}  \ \chi_L(t') \nabla \cdot s_0|v_a|^2 (t'')
\eeq

\noi so that
\bea
\label{2.20e}
&&\parallel \omega^{\alpha +n/2 \pm 0}  \ s_c \parallel \ \leq \ \int_t^1 \ dt'\ t{'}^{\gamma -2} \int_0^{t'} dt'' \parallel \omega^{\alpha + 2 + \gamma - n/2 \pm 0}  \ \chi_L(t') \ s_0|v_a|^2 (t'')\parallel \nn \\
&&\leq  C \ \int_t^1 dt' \  t{'}^{\lambda_{\alpha +1} - 2} \int_0^{t'}  dt'' \parallel \omega^{n/2 \pm 0} s_0(t'') \parallel \ \parallel v_a(t'') ; H^\rho \parallel^2
\eea

\noi by (\ref{2.3e}) (\ref{2.4e})
\bea
\label{2.21e}
&\leq&  C \ a_0^2 \ a_a^2 \int_t^1 dt'\ t{'}^{\lambda_{\alpha + 1}-2} \int_0^{t'} dt'' \ t{''}^{\lambda_0 - 1}  \nn \\
&\leq&  C \ a_0^2 \ a_a^2 \ \gamma^{ - 1} (1 - 2 \gamma )^{-1}  \  t^{\lambda_0 + \lambda_{\alpha + 1} -1}
\eea

\noi for $\alpha + 2 + \gamma >0$.\par \nobreak \hfill $\sq$ \par

In the applications, we shall use (\ref{2.12e}) (\ref{2.13e}) in the form
\beq
\label{2.22e}
\parallel \omega^{\alpha +n/2 \pm 0}  (s - s_0) \parallel \leq \ C \ a_0^2 \ a_a^2 \ t^{\lambda_0 +  \lambda_{\alpha + 1} -1}   
\eeq

\noi where we drop the dependence of the constant on $\gamma$. \par

In order to estimate the term $g_L(v) - g_L (v_a)$ in $L(v)$ (see (1.30)) we shall need estimates of $|v|^2 - |v_a|^2$. For that purpose, we shall use the fact that if $v$ satisfies the equation (\ref{1.22e}) for some real $V$ and if $v_a$ satisfies the equation (1.29), then the following formal conservation law holds
\beq
\partial_t \left ( |v|^2 - |v_a|^2\right ) = - \ {\rm Im}\ \overline{v} \Delta v + \nabla \cdot \left ( s |v|^2 - s_0 |v_a|^2 \right )
\label{2.23e}
\eeq

\noi (compare with (\ref{1.25e})  where $v_a$ satisfied (\ref{1.24e}) instead of (1.29)). We first give sufficient conditions for (\ref{2.23e}) to make sense and preliminary estimates which follow from it. The following lemma is a minor extension of Lemmas 2.6 and 2.7 of I.\\

\noi {\bf Lemma 2.6.} {\it Let $0 < \gamma < 1/2$ and $1 - \gamma /2 < \rho < n/2$. Let $0 < T \leq 1$, $I = (0, T]$ and $\overline{I} = [0, T]$. Let $v_0 \in H^\rho$. Let $v_a \in  {\cal C} (\overline{I}, H^\rho )$ be the solution of (1.29) with $v_a (0) = v_0$. Let $v \in L^\infty (I, H^\rho ) \cap {\cal C} (\overline{I}, L^2)$ satisfy the equation (\ref{1.22e}) in I for some real $V \in L_{loc}^\infty (I, L^\infty )$, with $v(0) = v_0$. Then 
\beq
|v(t)|^2 - |v_a(t)|^2 = V_1(t) + V_2(t)
\label{2.24e}
\eeq

\noi where
\bea 
\label{2.25e}
V_1 (t) &=& - \int_0^t dt' \ {\rm Im} \ \overline{v} \Delta v(t') \\
V_2 (t) &=& \nabla \cdot \int_0^t dt' \left ( s |v|^2 - s_0 |v_a|^2\right ) (t')
\label{2.26e}
\eea

\noi and for all $t \in I$, $V_1$, $V_2$ satisfy the estimates 
\beq
\label{2.27e}
\parallel \omega^{2\sigma - 2 - n/2}  \ V_1(t) \parallel \ \leq \ C \ a^2 \ t 
\eeq

\noi for $1/2 < \sigma \leq \rho \wedge (1 + n/4)$,
\beq  
\label{2.28e}
\parallel \omega^{2\sigma - 1 - n/2}  \ V_2(t) \parallel \ \leq \ C \ a^2 \ a_1^2 \ t^{\lambda_0} \left ( 1 + a^2 \  t^{\lambda_1} \right ) 
\eeq 

\noi for $0 < \sigma \leq \rho$, with}
$$a \ = \ \parallel v ; L^\infty (I, H^\rho ) \parallel \ , \quad a_1 = a \exp \left ( C\ a^2 \ T^{\lambda_1}\right ) \ . $$ 
\vskip 3 truemm

\noi {\bf Indication of proof.}

The proof is essentially the same as that of Lemma 2.7 of I. In particular the estimate (\ref{2.27e}) is identical with (\ref{2.22e}) of I. Here we give only the proof of the estimate (\ref{2.28e}), which is new. Note that $V_2$ here is more complicated than the corresponding $V_2$ of I. We estimate 
\bea  
\label{2.29e}
&&\parallel \omega^{2\sigma - 1 - n/2}  \ V_2(t) \parallel \ \leq \  \int_0^t dt' \parallel  \omega^{2\sigma - n/2} \left ( s|v|^2 - s_0 |v_a|^2 \right ) \parallel  (t') \nn \\
&&\leq \ C \int_0^t dt' \left ( \parallel \omega^{n/2 \pm 0}  \ s \parallel \ \parallel \omega^\sigma  v \parallel^2 \ + \ \parallel \omega^{n/2 \pm 0}  \ s_0 \parallel \ \parallel \omega^\sigma  v_a \parallel^2 \right ) (t')
\eea 

\noi by (\ref{2.3e}) (\ref{2.4e}), from which (\ref{2.28e}) follows by the use of (\ref{2.10e}) (\ref{2.22e}) and from the fact that $a_0 \leq a$ so that $a_a \leq a_1$. \par \nobreak \hfill $\sq$\par

The estimate (\ref{2.28e}) of $V_2$ is too rough for the subsequent applications. In particular it fails to exploit the expected cancellation between $s|v|^2$ and $s_0|v_a|^2$. In order to take the advantage of the latter, we rewrite 
\bea 
\label{2.30e}
V_2 &=& \nabla \cdot \int_0^t dt' \left ( (s - s_0)|v|^2\right ) (t') + \nabla \cdot \int_0^t dt' \left ( s_0 \left ( |v|^2 - |v_a|^2 \right ) \right ) (t')  \nn \\
&=& V_3 + V'_3
 \eea

\noi and we substitute again (\ref{2.24e}) in $V'_3$ so that $V'_3 = V_4 + V_5$ with 
 \bea 
\label{2.31e}
V_4 &=& \nabla \cdot \int_0^t dt'  s_0 (t') \ V_1(t') \nn \\
&=& - \nabla \cdot \int_0^t dt' \int_0^{t'} dt'' s_0(t') \ {\rm Im} \ \overline{v} \Delta v (t'') \ ,  \\
V_5 &=& \nabla \cdot \int_0^t dt'  s_0 (t') \ V_2(t') \nn \\
&=& \nabla \cdot \int_0^t dt' \int_0^{t'} dt'' s_0(t') \ \nabla \cdot  \left ( s |v|^2 - s_0 |v_a|^2 \right )  (t'') \ .  
\label{2.32e}
 \eea

\noi We then estimate $V_3$, $V_4$ and $V_5$ in the following lemma.\\

\noi {\bf Lemma 2.7.} {\it Let the assumptions of Lemma 2.6 be satisfied. Then the following estimates hold for all $t \in I$~: 
\beq  
\label{2.33e}
\parallel \omega^{2\sigma - 1 - n/2}  \ V_3(t) \parallel \ \leq \ C \ a^4 \ a_1^2 \ t^{\lambda_0 + \lambda_1} 
\eeq 

\noi for $0 < \sigma \leq \rho$, 
\bea  
\label{2.34e}
\left | <\psi , V_4 (t) > \right | &\leq& C \ a^4 \Big \{  t^{\lambda_0 + 1} \parallel \omega^{3 - 2\sigma + n/2} \ \psi \parallel \nn \\
&&+ \ \chi (\sigma \leq 1) \ t^{\lambda_\star + 1} \parallel \omega^{1 + n/2 \pm 0} \ \psi \parallel  \Big \}
\eea 

\noi for $1/2 < \sigma \leq \rho \wedge (1 + n/4)$ and for all $\psi$ such that the last two norms are finite whenever they occur, with 
\beq
\lambda_\star = \gamma - (1/2) (3 + \gamma - 2\sigma - 2 \rho )_+
\label{2.35e}
\eeq

\noi and $\chi (\sigma \leq 1) = 1$ (resp. {\rm 0}) if $\sigma \leq 1$ (resp. $\sigma > 1$).
\beq  
\label{2.36e}
\parallel \omega^{2\sigma - 2 - n/2}  \ V_5(t) \parallel \ \leq \ C \ a^4\ a_1^2\ t^{2 \lambda_0}  \left ( 1 + a^2 t^{\lambda_1}  \right )
\eeq 
 
\noi for $1/2 < \sigma \leq \rho$.}\\ 

\noi {\bf Proof.} We first estimate $V_3$. We obtain
\bea  
\label{2.37e}
&&\parallel \omega^{2\sigma - 1 - n/2}  \ V_3(t) \parallel \ \leq \ \int_0^t dt' \parallel \omega^{2\sigma - n/2} (s - s_0) |v|^2 \parallel (t') \nn \\ 
&&\leq \ C \int_0^t dt' \parallel \omega^{n/2 \pm 0} (s - s_0) \parallel \ \parallel \omega^{\sigma} v \parallel^2 (t')
\eea

\noi by (\ref{2.3e}) (\ref{2.4e}), from which (\ref{2.33e}) follows by the use of (\ref{2.22e}) and integration over time. \par

We next estimate $V_4$. Let $w =\ {\rm Im} \ \overline{v} \Delta v$. We rewrite
\beq  
\label{2.38e}
<\psi , V_4  > = \int_0^t dt' \int_0^{t'} dt'' \ <s_0(t') \cdot \nabla \psi , w (t'')>
\eeq 

\noi and we know (see the proof of (\ref{2.22e}) in Lemma 2.7 of I) that 
\beq  
\label{2.39e}
\parallel \omega^{2\sigma - 2 - n/2}  \ w \parallel \ \leq \ C \parallel \omega^{\sigma}  \ v \parallel^2
\eeq 

\noi for $1/2 < \sigma \leq \rho \wedge (1 + n/4)$, so that 
\beq  
\label{2.40e}
\left | <\psi , V_4  > \right |  \ \leq \ C \int_0^t dt' \int_0^{t'} dt'' \parallel \omega^{n/2 + 2 - 2\sigma} \ s_0(t') \nabla \psi \parallel \ \parallel \omega^{\sigma}  \ v (t'')\parallel^2 \ .
\eeq

 \noi We next estimate
\bea  
\label{2.41e}
&&\parallel \omega^{n/2+ 2 - 2\sigma} \ s_0(t') \nabla \psi \parallel \ \leq \ C \Big \{ \parallel \omega^{n/2 \pm 0} \ s_0(t') \parallel \ \parallel \omega^{n/2+ 3 - 2\sigma} \ \psi \parallel \nn \\
&&+ \chi (\sigma \leq 1) \parallel \omega^{n/2+ 2 - 2\sigma} \ s_0(t')  \parallel \ \parallel \omega^{1 + n/2 \pm 0} \psi \parallel  \Big \}
\eea 

\noi by (\ref{2.3e}) (\ref{2.4e}). The estimate (\ref{2.34e}) then follows from (\ref{2.40e}) (\ref{2.41e}) by the use of (\ref{2.10e}) and integration over time. The condition $3 - 2 \sigma + \gamma > 0$ needed to apply (\ref{2.10e}) in the second term of (\ref{2.41e}) is always fulfilled for $\sigma \leq 1$.

We finally estimate $V_5$. From (\ref{2.32e}) we obtain
\bea
\label{2.42e}
&&\parallel \omega^{2\sigma - 2 - n/2}  \ V_5 \parallel \ \leq \ \int_0^t dt' \parallel \omega^{2\sigma - 1 - n/2} \ s_0(t') \ V_2 (t') \parallel \nn \\ 
&&\leq \ C \int_0^t dt' \parallel \omega^{n/2 \pm 0} s_0(t') \parallel \ \parallel \omega^{2\sigma - 1 - n/2} \ V_2(t') \parallel
\eea

\noi for $1/2 < \sigma \leq \rho$ by (\ref{2.3e}), and (\ref{2.36e}) follows from (\ref{2.42e}) by the use of (\ref{2.10e}) (\ref{2.28e}) and integration over time.\par\nobreak \hfill $\sq$\par

For $\sigma > 1$, the estimate (\ref{2.34e}) of $V_4$ reduces to
\beq  
\label{2.43e}
\parallel \omega^{2\sigma - 3 - n/2}  \ V_4 \parallel \ \leq \ C \ a^4 \ t^{\lambda_0 + 1}  \ .
\eeq

\noi In the more interesting case $\sigma \leq 1$, it yields an estimate of $V_4$ in the space $\dot{H}^{2\sigma - 3 - n/2}  + \dot{H}^{-1 - n/2 \pm 0}$, but that space is not space dilation homogeneous. In the applications, we shall use (\ref{2.34e}) with a time dependent $\psi$, which will restore the space time dilation homogeneity of the estimate. More generally, we shall use the following lemma.\\

\noi {\bf Lemma 2.8.} {\it Let the assumptions of Lemma 2.6 be satisfied. \par

(1) Let $0 < \sigma ' < 1 + \gamma$ and let
\beq  
\label{2.44e}
\mu_j = \gamma - (1/2) \left ( j + 1 + \gamma - \sigma ' - 2 \rho \right )_+  \ .
\eeq

\noi Then the following estimates hold~:
\bea  
\label{2.45e}
&&t^{\gamma - 2} \parallel \omega^{\gamma - \sigma ' - n/2}  \ \chi_L \ V_1(t) \parallel \ \leq \ C \ a^2 \ t^{\mu_1 - 1}  \ , \\
\label{2.46e}
&&t^{\gamma - 2} \parallel \omega^{\gamma - \sigma ' - n/2}  \ \chi_L \ V_3(t) \parallel \ \leq \ C \ a^4 \ a_1^2 \ t^{2 \lambda_0 + \mu_1 - 2}  \ , \\
\label{2.47e}
&&t^{\gamma - 2} \parallel \omega^{\gamma - \sigma ' - n/2}  \ \chi_L \ V_4(t) \parallel \ \leq \ C \ a^4 \ t^{\lambda_0 + \mu_2 - 1}  \ , \\
&&t^{\gamma - 2} \parallel \omega^{\gamma - \sigma ' - n/2}  \ \chi_L \ V_5(t) \parallel \ \leq \ C \ a^4 \ a_1^2 \ t^{2 \lambda_0 + \mu_1 - 2} \left ( 1 + a^2 t^{\lambda_1} \right )  \ . 
\label{2.48e}
\eea

(2) Let $\sigma ' = 0$. Then the estimates (\ref{2.45e})-(\ref{2.48e}) hold with $\mu_j$ replaced by $\lambda_j$ in the right hand sides.} \\

\noi {\bf Proof} \underline{Part (1)}. We first estimate $V_1$. From (\ref{2.27e}) we obtain
\beq  
\label{2.49e}
t^{\gamma - 2} \parallel \omega^{\gamma - \sigma ' - n/2}  \ \chi_L \ V_1 \parallel \ \leq \ C \ a^2 \ t^{\gamma - 1 - (1/2) (2 + \gamma - \sigma ' - 2 \sigma )}  
\eeq

\noi under the conditions
\beq
\label{2.50e}
 \left\{ \begin{array}{l}
1 < 2\sigma \leq 2\rho \wedge (2 + n/2)  \\
\\
\sigma ' + 2 \sigma \leq 2 + \gamma
\end{array} \right .
\eeq

\noi (which make the condition $\sigma ' < 1 + \gamma$ unavoidable).\par

For $\sigma ' + 2\rho \geq 2 + \gamma$, we choose $\sigma$ so that $\sigma ' + 2\sigma  = 2 + \gamma$.\par

For $\sigma ' + 2\rho \leq 2 + \gamma$, we choose $\sigma = \rho$. That choice satisfies (\ref{2.50e}). In particular, the condition $2\sigma  \leq 2 + n/2$ follows from the condition $2\sigma \leq 2 + \gamma$. This proves (\ref{2.45e}).\par

We next estimate $V_3$. The estimate (\ref{2.33e}) is not sufficient for that purpose. We estimate instead for $\sigma ' < 2\sigma$
\begin{eqnarray*}
&&\parallel \omega^{2\sigma - \sigma '  - 1 - n/2}  \ V_3(t) \parallel \ \leq \ \int_0^t dt' \parallel \omega^{2\sigma - \sigma ' - n/2} \ (s - s_0) |v|^2  \parallel (t') \\ 
&&\leq \ C \int_0^t dt' \parallel \omega^{n/2 - \sigma '} (s - s_0) \parallel \ \parallel \omega^{\sigma} v \parallel^2 (t')
\end{eqnarray*}

\noi by (\ref{2.4e}),
\beq
\label{2.51e}
\leq \ C \ a^4 \ a_1^2 \ t^{\lambda_0 + \mu_1} 
\eeq

\noi by (\ref{2.22e}) (we do not need the $L^\infty$ norm of $s-s_0$, which allows us to use $\mu_j$ defined by (\ref{2.44e}) with $(\ )_+$ instead of $[\ ]_+$). From (\ref{2.51e}) we obtain
$$t^{\gamma - 2} \parallel \omega^{\gamma - \sigma ' - n/2}  \ \chi_L \ V_3 \parallel \ \leq \ C \ a^4 \ a_1^2 \ t^{\lambda_0 + \mu_1 + \gamma - (1/2) (1 + \gamma - 2\sigma ) - 2}$$

\noi for $0 < \sigma \leq \rho \wedge (1 + \gamma )/2$. Choosing $\sigma = \rho \wedge (1 + \gamma )/2$ yields (\ref{2.46e}).\par

We next estimate $V_4$. From (\ref{2.34e}) we obtain
\bea
\label{2.52e}
&&t^{\gamma - 2} \parallel \omega^{\gamma - \sigma ' - n/2}  \ \chi_L \ V_4 \parallel \ \leq \ C \ a^4 \ \mathrel{\mathop {\rm Sup}_{\parallel \psi \parallel = 1}} \ t^{\gamma - 1} \nn \\
&&\Big \{ t^{\lambda_0} \parallel \omega^{3 + \gamma - \sigma ' - 2\sigma}  \ \chi_L \ \psi \parallel \ + \ \chi (\sigma \leq 1) t^{\lambda_\star} \parallel \omega^{1 + \gamma - \sigma ' \pm 0} \ \chi_L \ \psi \parallel \Big \} \nn \\
&&\leq C\ a^4 \Big \{ t^{\lambda_0 - 1 + \gamma - (1/2) (3 + \gamma - \sigma ' - 2\sigma )} \nn \\
&&+ \chi (\sigma \leq 1) \ t^{\lambda_\star - 1 + \gamma - (1/2) (1 + \gamma - \sigma ')} \Big \}
\eea

\noi under the conditions
\beq
\label{2.53e}
 \left\{ \begin{array}{l}
1 < 2\sigma \leq 2\rho \wedge (2 + n/2)  \\
\\
\sigma ' + 2 \sigma \leq 3 + \gamma
\end{array} \right .
\eeq

\noi (together with the condition $\sigma ' < 1 + \gamma$). We next show that the second term in the last bracket in (\ref{2.52e}) is better behaved, namely has a larger time exponent than the first one. In fact 
\begin{eqnarray*}
&&2 \left ( \lambda_\star + \gamma - (1/2) \left ( 1 + \gamma - \sigma ') - \lambda_0 \right ) - 2 \gamma + (3 + \gamma - \sigma ' - 2 \sigma ) \right ) \\
&&= - (3 + \gamma - 2\sigma  - 2\rho )_+ + 2  - 2\sigma  + (1 + \gamma - 2\rho )_+ \geq 0
\end{eqnarray*}

\noi for $\sigma \leq 1$, since $(a + b)_+ \leq a_+ + b$ for $b \geq 0$. We can therefore omit that second term. We next choose $\sigma$. \par

For $\sigma ' + 2 \rho \geq 3 + \gamma$, we choose $\sigma$ so that $\sigma ' + 2\sigma = 3 + \gamma$. \par

For $\sigma ' + 2 \rho \leq 3 + \gamma$, we choose $\sigma = \rho$. That choice satisfies (\ref{2.53e}). In particular the condition $2\sigma > 1$ follows from $\sigma ' < 1 + \gamma$ in the first case and from $2\rho > 1$ in the second one. The condition $2\sigma \leq 2 + n/2$ follows from $2\sigma \leq 3 + \gamma$ for $n \geq 3$ and from $2\sigma  \leq  2 \rho < n$ for $n \leq 4$. This proves (\ref{2.47e}). \par

We finally estimate $V_5$. From (\ref{2.36e}) we obtain
\beq
\label{2.54e}
t^{\gamma - 2} \parallel \omega^{\gamma - \sigma ' - n/2}  \ \chi_L \ V_5 \parallel \ \leq \ C \ a^4 \ a_1^2  \left ( 1 + a^2\  t^{\lambda_1} \right )   t^{2\lambda_0 + \gamma - 2 - 1/2 (2 + \gamma - \sigma ' - 2\sigma )}
\eeq

\noi under the conditions
\beq
\label{2.55e}
 \left\{ \begin{array}{l}
1 < 2\sigma \leq 2\rho   \\
\\
\sigma ' + 2 \sigma \leq 2 + \gamma \ .
\end{array} \right .
\eeq

For $\sigma ' + 2 \rho \geq 2 + \gamma$, we choose $\sigma$ so that $\sigma ' + 2\sigma = 2 + \gamma$. \par

For $\sigma ' + 2 \rho \leq 2 + \gamma$, we choose $\sigma = \rho$. That choice satisfies (\ref{2.55e}). In particular the condition $2\sigma > 1$ follows from $\sigma ' < 1 + \gamma$ in the first case and from $2\rho > 1$ in the second one. This proves (\ref{2.48e}). \\

\noi \underline{Part (2)}. The proof is similar but simpler. \par\nobreak \hfill $\sq$ \par

\noi {\bf Remark 2.1.} The fact that $\sigma ' = 0$ in the proof of Part (2) and more generally the need of $L^\infty$ norms requires the use of $[\ ]_+$ in $\lambda_j$, whereas for $\sigma ' > 0$ one can use $(\ )_+$ in the definition of $\mu_j$.

\mysection{The linearized Cauchy problem for v} 
\hspace*{\parindent}
In this section we study the Cauchy problem for the linearized equation (1.35) with $L(v)$ defined by (1.30) for a given $v$, with initial time $t_0 \geq 0$. We first give a preliminary result with $t_0 > 0$, where we do not study the behaviour of the solution as $t$ tends to zero. \\

\noi {\bf Proposition 3.1.} {\it Let $0 < \gamma < 1$ and $\rho > 1 - \gamma /2$. Let $I = (0, T]$, let $v_0 \in H^\rho$ and let $v \in L_{loc}^\infty (I, H^\rho )$. Let $s$, $s_0$ and $v_a$ be defined by (1.31)-(1.34) and (1.29) with $v_a(0) = v_0$. Let $0 \leq \rho ' < n/2$, let $0 < t_0 \leq T$ and let $v'_0 \in H^{\rho '}$. Then the equation (1.35) has a unique solution $v' \in {\cal C} (I, H^{\rho '})$ with $v'(t_0) = v'_0$. The solution satisfies
$$\parallel v'(t) \parallel\ = \ \parallel v'_0 \parallel$$

\noi for all $t \in I$ and is unique in ${\cal C}(I, L^2)$.}\\

The proof is sketched in Appendix A2 of I.\par

We next study the boundedness and continuity properties near $t=0$ of the solutions of (1.35) obtained in Proposition 3.1. Since we shall eventually be interested in taking $\rho ' = \rho$, we already impose the condition $\rho < n/2$ in the next proposition (see however Remark 3.2 below).\\

\noi {\bf Proposition 3.2.} {\it Let $1/3 < \gamma < 1/2$ and $2 - 5\gamma /2 < \rho < n/2$. Let $I = (0, T]$ and $\overline{I} = [0, T]$, let $v_0 \in H^\rho$ and let $v \in L^\infty (I, H^\rho ) \cap {\cal C} (\overline{I}, L^2)$ with $v(0) = v_0$. Let $s$, $s_0$ and $v_a$ be defined by (1.31)-(1.34) and (1.29) with $v_a(0) = v_0$. Let $v$ satisfy the equation (\ref{1.22e}) in $I$ for some real $V \in L_{loc}^\infty (I, L^\infty )$. Let $1/2 \leq \rho ' < n/2$ and let $v' \in {\cal C} (I, H^{\rho '})$ be a solution of the equation (1.35) in $I$. Then \par

(1) $v' \in ({\cal C} \cap L^\infty ) (I, H^{\rho '}) \cap {\cal C}_w (\overline{I} , H^{\rho '}) \cap {\cal C} (\overline{I} , H^{\sigma})$ for $0 \leq \sigma < \rho '$.\par

(2) For all $t \in \overline{I}$, $t_1 \in I$, the following estimate holds 
\beq  
\label{3.1e}
\parallel \omega^{\rho '}  v'(t) \parallel\ \leq \ \parallel \omega^{\rho '}  v'(t_1) \parallel\ E(|t- t_1|)
 \eeq

\noi where
\bea  
\label{3.2e}
E(t) &\equiv& E(t,a) = \exp \left \{ C \ a^2 (1 + a^2) (1 + a_1^2)^2 \ t^{2\gamma + \lambda_1 - 1}  \right \} \ , \\
 a &=& \parallel v ; L^\infty (I, H^\rho ) \parallel\quad , \ a_1 = a \exp \left ( C\ a^2 \ T^{\lambda_1} \right )
\label{3.3e}
 \eea
\noi and $\lambda_j$ is defined by (\ref{2.9e}).\par

(3) For all $t$, $t_1 \in \overline{I}$, the following estimate holds} 
\beq  
\label{3.4e}
\parallel  v'(t) - v'(t_1)\parallel \ \leq C|t-t_1|^{\rho '\gamma  \wedge (3\gamma - 1)} (1 + a^2)^2 (1 + a_1^2)^2\parallel  v'(t_1) ; H^{\rho '}\parallel \ .
 \eeq

\noi {\bf Remark 3.1.} The estimate (\ref{3.1e}) for $t$, $t_1 \in I$ holds for $0 \leq \rho ' < n/2$, as will be clear from the proof. The condition $\rho ' \geq 1/2$ is needed to derive (\ref{3.4e}) which is used in turn to extend (\ref{3.1e}) to $t = 0$.\\

\noi {\bf Remark 3.2.} The assumption $\rho < n/2$ in Proposition 3.2 can be dispensed with at the expense of using slightly different estimates, which yield different powers of $t$ in (\ref{3.2e})  and (\ref{3.4e}).\\

\noi {\bf Proof.} We know already that the $L^2$- norm of $v'$ is conserved. The bulk of the proof consists in deriving the estimates (\ref{3.1e}) and (\ref{3.4e}) for $t$, $t_1 \in I$. We begin with (\ref{3.1e}). From (1.30) (1.35) we obtain 
\bea  
\label{3.5e}
\partial_t \parallel \omega^{\rho '} v' \parallel^2 &=&  {\rm Im} \ <v', [\omega^{2\rho '}, L(v)]v'>\nn \\
&=& \ {\rm Re} \ <v', [\omega^{2\rho '}, s] \cdot \nabla v'>\ + \ {\rm Im} \ <v', [\omega^{2\rho '}, f] v'>
 \eea

\noi where 
\beq  
\label{3.6e}
f = (1/2) \left ( |s|^2 - |s_0|^2 \right ) + t^{\gamma - 2} g_S (v) +  t^{\gamma - 2} \left ( g_L (v) - g_L(v_a) \right ) \ .
 \eeq

\noi We estimate the first term in the RHS of (\ref{3.5e}) by Lemma 2.4 with $\lambda = 2\rho '$, $\alpha_1 = 0$, $\alpha_2 = 1$, $r_i = 2$, $\sigma_1 = \sigma_2 = \rho '$ and $\nu = 1$, so that $\sigma_0 = 1 + n/2$ and $q_0 = \infty$. \par

We estimate similarly the last term by Lemma 2.4 with $\lambda =2 \rho '$, $\alpha_1 = \alpha_2 = 0$, $r_i = 2$, $\sigma_1 = \sigma_2 = \rho '$ and $\nu = 1$, so that $\sigma_0 = n/2$ and $\delta (q_0) = n/2-1$. We obtain
\beq  
\label{3.7e}
\left |\partial_t \parallel \omega^{\rho '} v' \parallel^2 \right | \ \leq \ C\left ( \parallel \omega^{n/2} \nabla s \parallel\ + \ \parallel \nabla s \parallel_\infty \ + \ \parallel \omega^{n/2} f \parallel \right ) \parallel \omega^{\rho '} v' \parallel^2 \ .
\eeq

\noi We estimate the various norms successively. We first estimate $\nabla s$ by (\ref{2.10e}) (\ref{2.22e}) with $\alpha = 1$ so that
\beq  
\label{3.8e}
\parallel \omega^{n/2} \nabla s \parallel\ + \ \parallel \nabla s \parallel_\infty \ \leq \ C\parallel \omega^{n/2+ 1 \pm 0}  s \parallel\ \leq \ C\ a^2\left \{  t^{\lambda_1 - 1} + a_1^2 \ t^{\lambda_0 + \lambda_2 - 1} \right \} \ .
 \eeq

\noi We next estimate the contribution of $f$. From (\ref{2.10e}) (\ref{2.22e}) with $\alpha = 0$, we obtain
\bea  
\label{3.9e}
&&\parallel \omega^{n/2} \left ( |s|^2 - |s_0|^2 \right )\parallel\  \leq \ C\parallel \omega^{n/2\pm 0} (s + s_0) \parallel \ \parallel \omega^{n/2\pm 0} (s - s_0) \parallel \nn \\
&&\leq \ C\ a^4 \ a_1^2 \ t^{2\lambda_0 + \lambda_1 - 2} \left ( 1 + a_1^2\ t^{\lambda_1} \right ) \ .
 \eea

We next estimate 
\bea  
\label{3.10e}
t^{\gamma - 2} \parallel \omega^{n/2} g_S(v) \parallel &\leq& C\ t^{\gamma - 2+ \rho - \gamma /2} \parallel \omega^\rho v \parallel^2\nn \\
&\leq &C\ a^2\ t^{\lambda_1 - 1} 
 \eea	

\noi for $\rho \geq \gamma /2$.\par 

The contribution of the last term in $f$ is estimated by the use of Lemma 2.8, part (2). We obtain
\bea  
\label{3.11e}
&&t^{\gamma - 2} \parallel \omega^{n/2} \left ( g_L(v) - g_L(v_a) \right ) \parallel \ = \ C\ t^{\gamma - 2} \parallel \omega^{\gamma - n/2}\ \chi_L \left ( V_1 + V_3 + V_4 + V_5\right )  \parallel\nn \\
&&\leq C\  a^2\left \{  t^{\lambda_1 - 1} + a^2\ t^{\lambda_0 + \lambda_2 - 1} + a^2\ a_1^2 \ t^{2\lambda_0 + \lambda_1 - 2} \left ( 1 + a^2\ t^{\lambda_1} \right )\right \} \ .
 \eea	

\noi Collecting  (\ref{3.7e})-(\ref{3.11e}), we obtain
\beq  
\label{3.12e}
\left | \partial_t \parallel \omega^{\rho '} v' (t) \parallel^2\right |  \ \leq \ N(t) \parallel \omega^{\rho '} v'(t) \parallel^2 
\eeq

\noi where
\beq  
\label{3.13e}
N(t) = C\ a^2\left \{ t^{\lambda_1 - 1} + a_1^2 \ t^{\lambda_0 + \lambda_2 - 1} + a^2\ a_1^2 \ t^{2\lambda_0 + \lambda_1 - 2} \left ( 1 + a_1^2\ t^{\lambda_1} \right )\right \} \ .
 \eeq	

\noi In order to estimate $\parallel \omega^{\rho '} v'(t)\parallel$, we need $N(t)$ to be integrable in time at $t = 0$. We first compare the various time exponents occurring in (\ref{3.13e}), assuming only that $0 < \gamma < 1/2$ and $\rho > 1 - \gamma /2$, which is equivalent to $\lambda_1 > 0$. Clearly
$$2\lambda_0 + \lambda_1 - 1 < \lambda_1 \wedge (2  \lambda_0 + 2 \lambda_1 - 1)\ .$$

\noi Moreover, from $[a + b]_+ \leq [a_+]_+ + b$ for $b > 0$, it follows that 
\beq
\label{3.14e}
[3 + \gamma - 2 \rho ]_+ \leq [2 + \gamma - 2 \rho ]_+ + 2 - 2 \gamma
\eeq

\noi and therefore $\lambda_0 + \lambda_2 \geq 2 \lambda_0 + \lambda_1 - 1$. Keeping in (\ref{3.13e}) the dominant power of $t$ and using the fact that $\lambda_0 = \gamma$ under the previous assumptions, we obtain 
\beq  
\label{3.15e}
N(t) \leq C\  a^2 (1 + a^2) (1 + a_1^2)^2 \ t^{2\gamma + \lambda_1 - 2} \ .
 \eeq	

\noi The integrability condition of $N(t)$ at $t=0$ then becomes $2\gamma + \lambda_1 -1 > 0$ or equivalently $\gamma > 1/3$ and $\rho > 2 - 5\gamma /2$. \par

The estimate (\ref{3.1e}) (\ref{3.2e}) for $t_1$, $t \in I$ follows from (\ref{3.12e}) (\ref{3.15e}) by integration.\par

We next derive the estimate (\ref{3.4e}) for $t$, $t_1 \in I$. For that purpose we define (see  (\ref{1.8e}))
\beq  
\label{3.16e}
\widetilde{v}'(t) = U(-t) v'(t)
 \eeq
\beq  
\label{3.17e}
\widetilde{L} = L(v) + (1/2) \Delta = is \cdot \nabla + (i/2) (\nabla \cdot s ) + f
 \eeq

\noi with $f$ given by (\ref{3.6e}). We rewrite (1.35) as
\beq  
\label{3.18e}
i \partial_t \widetilde{v}'  = U(-t) \widetilde{L} U(t) \widetilde{v}'
 \eeq	

\noi so that for $t$, $t_1 \in I$, for fixed $t_1$,
\bea
\label{3.19e}
\partial_t \parallel \widetilde{v}' (t) - \widetilde{v}' (t_1) \parallel^2 &=&2 \ {\rm Im} \ < \widetilde{v}' (t) - \widetilde{v}' (t_1) , U(-t) \widetilde{L} \ U(t) \  \widetilde{v}' (t_1) >\nn \\
&=&2 \ {\rm Im} \ < w, \widetilde{L}\  v_\star>
\eea

\noi where
\beq
\label{3.20e}
\left \{ \begin{array}{l} v_\star (t) = U(t-t_1) v'(t_1) \\ \\ w(t) = v'(t) - v_\star (t) \ . \end{array} \right . 
\eeq

\noi We estimate 
\bea  
\label{3.21e}
&&\left |\partial_t \parallel w \parallel^2 \right | \ \leq \ 2\left | {\rm Re} \ <w, s\cdot \nabla v_*> \right |\nn \\
&&+\ C  \parallel w\parallel \left ( \parallel \omega^{n/2 - \sigma '} \nabla \cdot s \parallel\ + \ \parallel \omega^{n/2 - \sigma '} f \parallel \right )  \parallel \omega^{\sigma '} v'(t_1) \parallel 
\eea

\noi for some $\sigma '$ with $0 < \sigma ' \leq \rho '$, to be chosen later. \par

For $0 < \rho ' < 1$, we write
\beq
\label{3.22e}
<w,s\cdot \nabla v_*> \ = \ - <\omega^{- \rho '} \nabla \cdot sw, \omega^{\rho '} v_*>
\eeq

\noi and we estimate by Lemma 2.2
\beq
\label{3.23e}
\left |< w, s\cdot \nabla v_*>\right |  \ \leq \ C \parallel \omega^{1 - \rho '} w \parallel\left (  \parallel s \parallel_\infty \ + \ \parallel \omega^{n/2} s \parallel \right ) \ \parallel \omega^{\rho '} v'(t_1) \parallel \ .
\eeq

\noi For $\rho ' = 1$, we estimate
\beq
\label{3.24e}
\left |< w, s\cdot \nabla v_*>\right |  \ \leq \  \parallel w \parallel\ \parallel s \parallel_\infty \  \parallel \omega^{\rho '} v'(t_1) \parallel  \ .
\eeq

\noi For $\rho ' > 1$, we estimate
\beq
\label{3.25e}
\left |< w, s\cdot \nabla v_*>\right |  \ \leq \ C \parallel w \parallel\ \parallel \omega^{n/2- \sigma '} \nabla s \parallel  \ \parallel \omega^{\sigma '} v'(t_1) \parallel 
\eeq

\noi for $1 < \sigma ' \leq \rho '$.\par

Collecting (\ref{3.21e})-(\ref{3.25e}) yields 
\bea  
\label{3.26e}
&&\left |\partial_t \parallel w \parallel^2 \right | \ \leq \ C\Big \{ \chi  (\rho ' \leq 1) \parallel \omega^{1 - \rho '} w \parallel\left (  \parallel s \parallel_\infty \ + \ \parallel \omega^{n/2} s \parallel\right ) \parallel \omega^{\rho '} v'(t_1) \parallel \nn \\
&&+\   \parallel w\parallel \left ( \parallel \omega^{n/2- \sigma '} \nabla s \parallel\ + \ \parallel \omega^{n/2- \sigma '} f \parallel \right )  \parallel \omega^{\sigma '} v'(t_1) \parallel \Big \}
\eea

\noi with $0 < \sigma ' \leq \rho '$ and $\sigma ' > 1$ in the $\nabla s$ term if $\rho ' > 1$.\par

For $1/2 \leq \rho ' \leq 1$, we interpolate 
$$\parallel \omega^{1 - \rho '} w \parallel \ \leq \ y^\theta \parallel \omega^{\rho '} w \parallel^{1/\rho ' - 1}$$

\noi where
$$y = \ \parallel w \parallel^2 \qquad , \quad \theta = 1 - 1/(2 \rho ')$$

\noi so that (\ref{3.26e}) becomes
\bea  
\label{3.27e}
\left |\partial_t y \right | &\leq& \ C\left \{ \chi  (\rho ' \leq 1) \left ( \parallel s \parallel_\infty \ + \ \parallel \omega^{n/2} s \parallel\right ) a{'}_1^{1/\rho '} y^\theta \right . \nn \\
&&\left . +\   \left ( \parallel \omega^{n/2- \sigma '} \nabla s \parallel \ + \ \parallel \omega^{n/2- \sigma '} f \parallel \right )  a'_1 \ y^{1/2}  \right \} 
\eea

\noi with $a'_1 = \parallel v'(t_1);H^{\rho '}\parallel$. We estimate the various norms in (\ref{3.27e}) successively. From (1.31) (\ref{2.10e}) (\ref{2.22e}) with $\alpha = 0$ we obtain
\beq
\label{3.28e}
\parallel s\parallel_\infty \ + \ \parallel \omega^{n/2} s \parallel \ \leq\ C\ a^2 \   t^{\lambda_0 - 1} \left ( 1  + a_1^2 \ t^{\lambda_1} \right )
\eeq

\noi and with $\alpha = 1 - \sigma '$, we obtain
\beq
\label{3.29e}
\parallel \omega^{n/2- \sigma '} \nabla s \parallel \ \leq \ C\ a^2 \left ( t^{\mu_1- 1} + a_1^2 \ t^{\lambda_0 + \mu_2 - 1} \right )
\eeq

\noi for $0 < \sigma ' < 2 + \gamma$ (see Remark 2.1).\par

\noi We next estimate the contribution of $f$. From (\ref{2.10e}) (\ref{2.22e}) we obtain 
\bea
\label{3.30e}
&&\parallel \omega^{n/2- \sigma '} \left ( |s|^2 - |s_0|^2 \right ) \parallel \ \leq \ C\parallel \omega^{n/2 \pm 0} (s + s_0) \parallel \
\parallel \omega^{n/2 - \sigma '} (s - s_0) \parallel \nn \\
&&\leq \ C\ a^4\ a_1^2\ t^{2\lambda_0 + \mu_1 - 2} \left ( 1 + a_1^2 \  t^{\lambda_1} \right ) 
\eea

\noi for $0 < \sigma ' < (2 + \gamma ) \wedge n$. We then estimate
\beq
\label{3.31e}
t^{\gamma -2}  \parallel \omega^{n/2 - \sigma '} g_S(v) \parallel \ \ \leq \ C\ a^2\  t^{\mu_1 -1}\ , 
\eeq

\noi for $0 < \sigma ' \leq 2 \rho - \gamma$.\par

The contribution of the last term in $f$ is estimated by the use of Lemma 2.8, part (1) for $0 < \sigma ' < 1 + \gamma$ and $\sigma ' \leq \rho '$~:
\bea
\label{3.32e}
&&t^{\gamma -2}  \parallel \omega^{n/2 - \sigma '} \left ( g_L(v) - g_L(v_a)\right )  \parallel \ = \  C\  t^{\gamma -2} \parallel \omega^{\gamma - n/2- \sigma '} \chi_L \left ( V_1 + V_3 + V_4 + V_5 \right ) \parallel \nn \\
&&\leq \ C\ a^2\ \left \{ t^{\mu_1 - 1} + a^2\ t^{\lambda_0 + \mu_2 - 1} + a^2\ a_1^2\ t^{2\lambda_0 + \mu_1 - 2}\left ( 1 + a^2\  t^{\lambda_1}\right ) \right \} \ .
\eea

\noi Collecting (\ref{3.27e})-(\ref{3.32e}), we obtain
\bea  
\label{3.33e}
&&\left |\partial_t y \right | \ \leq  \ C\Big \{ \chi  (\rho ' \leq 1)  a^2 \ t^{\lambda_0 - 1} \left ( 1 + a_1^2 \ t^{\lambda_1}\right ) a{'}_1^{1/\rho '} \ y^\theta  \nn \\
&&+  a^2  \left ( t^{\mu_1 -1} + a_1^2 \   t^{\lambda_0 + \mu_2 - 1} + a^2 \ a_1^2\ t^{2\lambda_0 + \mu_1 - 2} \left ( 1 + a_1^2\  t^{\lambda_1}\right ) \right )   a'_1\ y^{1/2} \Big \} \ .
\eea

\noi We next choose $\sigma '$ as large as possible, namely $\sigma ' = \rho ' \wedge (1 + \gamma - 0)$ and we compare the various time exponents that occur in the last term of (\ref{3.33e}), assuming only that $0 < \gamma < 1/2$ and $\rho > 1 - \gamma /2$, which is equivalent to $\lambda_1 > 0$. Clearly $2\lambda_0 + \mu_1 - 1 < \mu_1 \wedge (2 \lambda_0 + \lambda_1 + \mu_1 - 1)$. We next show that $\lambda_0 + \mu_2 - 1 \geq 2 \lambda_0 + \mu_1 - 2$, or equivalently 
\beq  
\label{3.34e}
\left ( 3 + \gamma - 2 \rho - \rho '  \wedge (1 + \gamma - 0)\right )_+ \leq \left ( 2 + \gamma -2\rho - \rho ' \right )_+ + 2(1 - \gamma )\ ,
\eeq

\noi where we have used the fact that the limitation $\sigma ' < 1 + \gamma$ is not seen in $\mu_1$ for $\rho > 1/2$. The inequality (\ref{3.34e}) with $\rho '$ is proved in the same way as (\ref{3.14e}) and the inequality with $1 + \gamma$ is trivial since
$$2 - 2\rho < \gamma \leq 2(1 - \gamma )\ .$$

\noi The dominant exponent is then 
\begin{eqnarray*}  
2\lambda_0 + \mu_1 - 1 &=& 3\gamma - 1 - (1/2) \left ( 2 + \gamma - 2\rho - \rho' \right )_+ \\
&=& (3\gamma - 1) \wedge \left ( \rho - 2 + 5\gamma /2 +  \rho ' /2 \right ) \\
&\geq& \mu \equiv (3\gamma - 1) \wedge \rho ' /2 
\end{eqnarray*}

\noi for $\rho > 2 - 5 \gamma /2$. \par

From (\ref{3.33e}) we then obtain
\bea  
\label{3.35e}
\left |\partial_t y \right | &\leq&  C\Big \{ \chi  (\rho ' \leq 1)  a^2 \left ( 1 + a_1^2 \right ) t^{\gamma - 1} \ a{'}_1^{1/\rho '} \ y^\theta  \nn \\
&&+  a^2  (1 + a^2) \left ( 1 + a_1^2 \right )^2 t^{\mu -1} \   a'_1\ y^{1/2} \Big \} \ .
\eea

\noi Using the fact that the differential inequality 
$$\left |\partial_t y \right | \leq  \sum_i  \ b_i\ t^{\nu_i -1} \ y^{\theta_i}$$

\noi with $0 \leq \theta_i < 1$ and $\nu_i > 0$ implies
$$y(t) \leq  C \ \sum_i  \left (  b_i\ \nu_i^{-1} \left | t^{\nu_i} - t^{\nu_i}_1\right | \right )^{1/(1 - \theta_i)}$$

\noi for $t$, $t_1 > 0$ and $y(t_1) = 0$, we obtain
\beq  
\label{3.36e}
y(t)\leq  C\left \{ \chi  (\rho ' \leq 1)  \left ( a^2(1 + a_1^2)\right )^{2\rho '} \left | t^\gamma - t_1^\gamma\right |^{2\rho '} + \left ( a^2(1 + a^2) (1 + a_1^2)^2 \right )^2 \left | t^\mu - t_1^{\mu} \right |^2\right \}  a{'}_1^2
\eeq

\noi so that
\beq  
\label{3.37e}
\parallel w(t)\parallel \leq  C (1 + a^2)^2 (1 + a_1^2)^2 \left \{ \chi  (\rho ' \leq 1)  | t - t_1|^{\gamma\rho '} + |t - t_1|^\mu \right \}  a{'}_1 \ .
\eeq

\noi On the other hand
\bea  
\label{3.38e}
&&\parallel v'(t) - v'(t_1) \parallel \ \leq \ \parallel w(t) \parallel \ + \ \parallel \left ( U(t-t_1) - 1\right ) v'(t_1) \parallel \nn \\
&&\leq \ \parallel w(t) \parallel\ + \ |t-t_1|^{(\rho '/2) \wedge 1} \parallel \omega^{\rho ' \wedge 2} \ v'(t_1) \parallel \ .
\eea

\noi Collecting (\ref{3.37e}) (\ref{3.38e}) yields (\ref{3.4e}) for $t$, $t_1 \in I$. \par

We now exploit (\ref{3.1e}) and (\ref{3.4e}) in $I$ to complete the proof of the proposition. From (\ref{3.1e}) it follows that $v' \in L^\infty (I, H^{\rho '})$. From (\ref{3.1e}) and (\ref{3.4e}) it then follows that $v'$ has a limit $v'(0)$ in $L^2$ and that (\ref{3.4e}) holds for $t$, $t_1 \in \overline{I}$. It then follows by a standard abstract argument that $v'(0) \in H^{\rho '}$, that $v' \in {\cal C}_w (\overline{I}, H^{\rho '}) \cap {\cal C} (\overline{I}, H^{\sigma})$ for $0 \leq \sigma < \rho '$, and that (\ref{3.1e}) holds for all $t \in \overline{I}$, $t_1 \in I$. \par\nobreak \hfill $\sq$\par

\noi {\bf Remark 3.3.} The integrability of $N(t)$ at $t =0$ requires stronger conditions on $\rho$ than the subcriticality condition $\rho > 1 - \gamma /2$, or equivalently $\lambda_1 > 0$. That condition suffices to control the contributions of $\nabla s_0$, of $g_S(v)$ and of $V_1$ to (\ref{3.7e}). The terms $\nabla (s - s_0)$ and $V_4$ require $\lambda_0 + \lambda_2 > 0$, or equivalently $\rho > 3/2 - 3\gamma /2 = 1 - \gamma /2 + (1/2 - \gamma )$. The worst terms are $|s|^2 - |s_0|^2$, $V_3$ and $V_5$ which require $2\lambda_0 + \lambda_1 > 1$, or equivalently $\gamma > 1/3$ and $\rho > 2 - 5\gamma /2 = 1 - \gamma /2 + (1 - 2\gamma )$. \\

We have not proved so far that $v' \in {\cal C} (\overline{I}, H^{\rho '})$. This is true but requires a separate argument. \\

\noi {\bf Proposition 3.3.} {\it Under the assumptions of Proposition 3.2, $v' \in {\cal C} (\overline{I}, H^{\rho '})$ and (\ref{3.1e}) holds for all $t$, $t_1 \in \overline{I}$.}\\

\noi The proof is identical with that of Proposition 3.3 of \cite{4r}. \par

We can now state the main result on the Cauchy problem for the linearized equation (1.35). \\

\noi {\bf Proposition 3.4.} {\it Let $1/3 < \gamma < 1/2$ and $2 - 5\gamma /2 < \rho < n/2$. Let $I = (0, T]$ and $\overline{I} = [0, T]$ and let $v \in L^\infty (I, H^\rho ) \cap  {\cal C} (\overline{I}, L^2)$ with $v(0) = v_0$. Let $s$, $s_0$ and $v_a$ be defined by (1.31)-(1.34) and (1.29) 
 with $v_a (0) = v_0$. Let $v$ satisfy the equation (\ref{1.22e}) in $I$ for some real $V \in L_{loc}^\infty (I, L^\infty )$. Let $1/2 < \rho ' < n/2$ and let $v'_0 \in H^{\rho '}$. Let $t_0 \in \overline{I}$. Then there exists a unique solution $v' \in  {\cal C} (\overline{I}, H^{\rho '})$ of the equation (1.35) with $v'(0) = v'_0$. Furthermore $v'$ satisfies the estimates (\ref{3.1e}) and (\ref{3.4e}) for all $t$, $t_1 \in \overline{I}$. The solution is actually unique in ${\cal C} (\overline{I}, L^2)$.}\\

 The proof is identical with that of Proposition 3.4 of I.

\mysection{The nonlinear Cauchy problem at time zero for v and u$_{\bf c}$}
\hspace*{\parindent} 
In this section we prove that the nonlinear equation (\ref{1.19e}) for $v$, with $L(v)$ defined by (1.30), with initial data at time zero, has a unique solution in a small time interval. We then rewrite that result in terms of $u_c$, related to $v$ by (\ref{1.12e}), and we give some additional bounds and regularity properties for $u_c$. In order to solve the equation (\ref{1.19e}) for $v$, we show that the map $\Gamma : v \to v'$ defined by Proposition 3.4 with $t_0 = 0$ is a contraction. For that purpose, we need to estimate the difference of two solutions of the linearized equation (1.35). For any pair of functions or operators $(f_1, f_2)$, we define
$$f_\pm = (1/2) \left ( f_2 \pm f_1 \right ) \ .$$
\vskip 3 truemm

\noi {\bf Lemma 4.1.} {\it Let $1/3 < \gamma < 1/2$ and $2 - 5\gamma /2 < \rho < n/2$. Let $I = (0,T]$ and let $v_i$, $i = 1,2$ satisfy the assumptions of Proposition 3.4 with $v_i(0) = v_0$. Let $1/2 < \rho ' < n/2$ and let $v'_i$, $i =1,2$ be the solutions of the equation (1.35) with $v'_i (0) = v'_0 \in H^{\rho '}$ obtained in Proposition 3.4. Then the following estimate holds for all $t \in I$~:
 \beq
\label{4.1e}
\parallel v'_-; L^\infty ((0, t], H^{\rho '}) \parallel \ \leq \ C\ E(t,a) a(1 + a^2)^2 \left ( 1 + a_1^2\right ) a' \ t^{2\gamma +\lambda_1 - 1}  \parallel v_-; L^\infty ((0, t], H^\rho ) \parallel
\eeq 

\noi where $E(t, a)$ is defined by (\ref{3.2e}) and}
 \beq
\label{4.2e} 
a = \ {\rm Max}\parallel v_i; L^\infty (I, H^\rho ) \parallel\ , \quad a' = \ {\rm Max}\parallel v'_i; L^\infty (I, H^{\rho '} ) \parallel\  .
\eeq

\noi {\bf Proof.} From (1.35) we obtain
$$i\partial_t v'_- = L_2 \ v'_- + L_- \ v'_1$$

\noi where $L_i =L(v_i)$, $g_i = g(v_i)$, so that
$$L_- = t^{\gamma - 2}\ g_-  \ .$$

\noi We estimate for $0 \leq \sigma  \leq \rho '$
\beq
\label{4.3e}
\partial_t \parallel \omega^{\sigma } v'_- \parallel^2 \ = \ 2\ {\rm Im}  \left (  < \omega^{\sigma }  v'_-, \omega^{\sigma } L_2 v'_- > \ + \ <   \omega^{\sigma } v'_-, \omega^{\sigma } L_- v'_1 > \right )\ .
\eeq 

\noi By the estimates in the proof of Proposition 3.2 (see in particular  (\ref{3.1e}) ; see also Remark 3.1), we obtain
\beq
\label{4.4e}
\parallel \omega^{\sigma } v'_- (t)\parallel \ \leq \ E(t,a) \int_0^t dt'\ t{'}^{\gamma -2} \parallel  \omega^{\sigma } g_- \ v'_1(t')\parallel \ .
\eeq

\noi We next estimate
\bea
\label{4.5e}
&&\parallel \omega^{\sigma } g_- \ v'_1 \parallel \ \leq \ C \parallel \omega^{n/2 \pm 0} g_- \parallel \ \parallel \omega^\sigma  v'_1\parallel  \ , \\
&&\parallel \omega^{n/2 \pm 0} g_{S-}  \parallel \ \leq \ C \ t^{\rho - \gamma /2} \parallel  \omega^{\rho} v_-\parallel \ \parallel  \omega^{\rho} v_+\parallel\ .
\label{4.6e}
\eea 

\noi In order to estimate $g_{L-}$, we use again Lemma 2.8, part (2). From the conservation law (\ref{2.23e}) and from the fact that $v_- (0) = 0$ we obtain (see (\ref{2.24e}))
\beq
\label{4.7e}
|v(t)|_-^2 \equiv (1/2) \left ( |v_2|^2 - |v_1|^2 \right ) = V_{1-}(t) + V_{2-}(t)
\eeq

\noi where (see (\ref{2.25e}) (\ref{2.26e}))
\begin{eqnarray*}
&&V_{1-}(t) = - \int_0^t dt'\ {\rm Im} \left ( \overline{v}_+ \Delta v_- + \overline{v}_- \Delta v_+\right ) (t') \ , \\
&& V_{2-}(t) = \int_0^t dt'\ \nabla \cdot s|v|^2_-  (t') \ .
\end{eqnarray*}

\noi In the same way as in Section 2 (see (\ref{2.30e})-(\ref{2.32e})), we rewrite
\beq
\label{4.8e}
V_{2-} = \nabla \cdot \int_0^t dt' (s - s_0) |v|_-^2(t') + \nabla \cdot \int_0^t dt' s_0 |v|_-^2(t') 
 \eeq

\noi and we substitute again (\ref{4.7e}) into the last term, so that
\beq
V_{2_-} = V_{3_-} + V_{4_-} + V_{5_-}
\label{4.9e}
\eeq

\noi with
\bea
\label{4.10e}
V_{3_-} &=& \nabla \cdot \int_0^t dt' (s - s_0) |v|_-^2 (t') \ , \\
\label{4.11e}
V_{4_-} &=& \nabla \cdot \int_0^t dt' s_0(t') V_{1_-} (t') \nn \\
&=& - \nabla \cdot \int_0^t dt' \int_0^{t'} dt'' s_0 (t')  \ {\rm Im} \left ( \overline{v}_+ \Delta v_- + \overline{v}_- \Delta v_+\right ) (t'') \ , \\
V_{5_-} &=& \nabla \cdot \int_0^t dt' s_0(t') V_{2_-} (t') = \nabla \cdot \int_0^t dt' \int_0^{t'} dt'' s_0 (t') \nabla \cdot s|v|^2_- (t'') \ .
\label{4.12e}
\eea

\noi By essentially the same estimates as in Lemma 2.8, part (2), we obtain from (\ref{4.5e}) (\ref{4.6e})
\bea
\label{4.13e}
&&t^{\gamma - 2} \parallel \omega^{\sigma } g_- \ v'_1 \parallel \ \leq \ C \Big \{ a \  t^{\lambda_1 - 1} + a^3 \ t^{\lambda_0 + \lambda_2 - 1 }\nn \\
&&+ a^3\ a_1^2 \  t^{2\lambda_0 + \lambda_1 - 2} \left ( 1 + a^2 \ t^{\lambda_1} \right ) \Big \} a'\parallel v_-; L^\infty ((0, t], H^{\rho} ) \parallel \nn \\
&&\leq C\ a(1 + a^2)^2 (1 + a_1^2) a'\ t^{2 \gamma + \lambda_1 - 2} \parallel v_-; L^\infty ((0, t], H^{\rho} ) \parallel
\eea

\noi for $0 \leq \sigma \leq \rho '$, by keeping the dominant power of $t$ in the last inequality (see (\ref{3.13e})-(\ref{3.15e})). \par 

Substituting (\ref{4.13e}) into (\ref{4.4e}) yields (\ref{4.1e}). \par\nobreak \hfill $\sq$ \par

We can now state the main result on the Cauchy problem at time zero for the equation (\ref{1.19e}). \\

\noi {\bf Proposition 4.1.} {\it Let $1/3 < \gamma  < 1/2$ and $2 - 5\gamma /2 < \rho < n/2$, let $v_0 \in H^\rho$ and let $s$, $s_0$ and $v_a$ be defined by (1.31)-(1.34) and (1.29) with $v_a(0) = v_0$. Then there exists $T > 0$ and there exists a unique solution $v \in {\cal C}([0, T], H^\rho )$ of the equation (\ref{1.19e}) with $L(v)$ defined by (1.30), with $v(0) = v_0$. One can ensure that 
  \bea
\label{4.14e}
&& \parallel v; L^\infty ([0, T], H^\rho )  \parallel \ \leq \ R=2 \parallel v_0; H^\rho \parallel \\
&&C\ R^2 \left ( 1 + R^2 \right )^3 T^{2\gamma + \lambda_1 - 1} = 1
\label{4.15e}
\eea

\noi for some $C$ independent of $v_0$.}\\

The proof is identical with that of Proposition 4.1 of I.\par

We finally translate the main result of Proposition 4.1 in terms of $u_c$ and we derive additional bounds and regularity properties for $u_c$.\\

\noi {\bf Proposition 4.2.} {\it Let $1/3 < \gamma < 1/2$ and $2 - 5\gamma /2 < \rho < n/2$, let $v_0 \in H^\rho$. Let $\varphi$ be defined by (1.28) with $\varphi (1) = 0$, with $s_0$ defined by (1.32) and $v_a$ defined by (1.29) with $v_a(0) = 0$. Then there exists $T > 0$ and there exists a unique solution $u_c \in {\cal C}((0, T], H^\rho )$ of the equation (\ref{1.11e}) such that $v$ defined by (\ref{1.12e}) satisfies the equation (\ref{1.19e}) with $L(v)$ defined by (1.30), with $v(0) = v_0$. Furthermore $u_c$ satisfies the estimate 
 \beq
\label{4.16e}
\parallel u_c(t); H^\rho  \parallel \ \leq  C\ a_0 \left ( 1 + a_0^2 \left ( 1 + a_0^2\ t^\gamma\right )  t^{\gamma - 1} \right )^{1 + [\rho ]} 
\eeq

\noi for all $t \in (0, T]$, where $[ \rho ]$ is the integral part of $\rho$ and}
$$a_0 \ = \ \parallel v_0; H^\rho  \parallel \ .$$
\vskip 3 truemm

\noi {\bf Sketch of proof.} The proof is the same as that of Proposition 4.2 of I, using the appropriate estimates of $\varphi$. As in the latter we estimate
 \beq
\label{4.17e}
\parallel \omega^{\rho } u_c \parallel \ \leq \ C \left ( 1 \ + \ \parallel \omega^{n/2} \varphi \parallel\right )^{1+[\rho ]} \parallel \omega^{\rho } v\parallel \ . 
\eeq 

\noi In the present case, $\varphi$ can be written as $\varphi = \varphi_0 + \varphi_b + \varphi_c$ where the various terms are defined in analogy with (1.32)-(1.34) with $\nabla$ omitted and are estimated in the same way as in Lemma 2.5, parts (1) and (3). One obtains 
\beq
\label{4.18e}
\parallel \omega^{n/2} \varphi \parallel \ \leq \ C \left ( a_0^2 \ t^{\gamma - 1} + a_0^2 \ a_1^2 \ t^{2\gamma - 1}   \right ) \eeq
 
\noi which together with (\ref{4.17e}) and with the fact that $a_1 \leq C\ a_0$ under the condition (\ref{4.15e}), implies (\ref{4.16e}).\par\nobreak \hfill $\sq$ \par

\end{document}